\documentclass[aap,seceqn,nameyear,MSNbibl,dvips]{arximspdf}

\doi{10.1214/10-AAP712}
\volume{21}
\issue{2}
\pubyear{2011}
\firstpage{464}
\lastpage{483}

\makeatletter

\newcommand{\rr}{\mathbb{R}}
\newcommand{\hatK}{\hat{K}}
\newcommand{\la}{\lambda}

\newtheorem{theorem}{Theorem}[section]
\newtheorem{lemma}{Lemma}[section]
\newproclaim{remark}{Remark}[section]

\makeatother

\begin{document}
\begin{frontmatter}

\title{Nonnormal approximation by
Stein's method of exchangeable pairs with application to
the Curie--Weiss model}
\runtitle{Nonnormal approximation}

\begin{aug}
\author[A]{\fnms{Sourav} \snm{Chatterjee}\thanksref{t1}\ead
[label=e1]{sourav@cims.nyu.edu}} and
\author[B]{\fnms{Qi-Man} \snm{Shao}\corref{}\thanksref{t2}\ead
[label=e2]{maqmshao@ust.hk}}
\runauthor{S. Chatterjee and Q.-M. Shao}
\affiliation{New York University and Hong Kong
University of Science and Technology}
\address[A]{Department of Mathematics\\
Courant Institute of Mathematical Sciences\\
New York University\\
251 Mercer Street, New York\\
New York 10012\\
USA\\
\printead{e1}}
\address[B]{Department of Mathematics\\
Hong Kong University\\
of Science and Technology\\
Clear Water Bay, Kowloon, Hong Kong\\
China\\
\printead{e2}}
\end{aug}

\thankstext{t1}{Supported by NSF Grant DMS-07-07054 and a Sloan
Research Fellowship.}

\thankstext{t2}{Supported in part by Hong Kong RGC 602206, 602608 and RGC 603710.}

\received{\smonth{12} \syear{2008}}
\revised{\smonth{4} \syear{2010}}

%
\begin{abstract}
Let $(W, W')$ be an exchangeable pair. Assume that
\[
E(W-W' | W) = g(W)+r(W),
\]
where $g(W)$ is a dominated term and $r(W)$ is negligible. Let $G(t) =
\int_0^t g(s) \,ds $ and define $p(t) = c_1 e^{-c_0 G(t)}$, where $c_0$
is a properly chosen constant and $c_1 = 1/\int_{-\infty}^\infty
e^{-c_0 G(t)} \,dt$. Let $Y$ be a random variable with the probability
density function $p$. It is proved that $W$ converges to $Y$ in
distribution when the conditional second moment of $(W-W')$ given $W$
satisfies a law of large numbers. A Berry--Esseen type bound is also
given. We use this technique to obtain a Berry--Esseen error bound of
order $1/\sqrt{n}$ in the noncentral limit theorem for the
magnetization in the Curie--Weiss ferromagnet at the critical
temperature. Exponential approximation with application to the spectrum
of the Bernoulli--Laplace Markov chain is also discussed.
\end{abstract}

%

\end{frontmatter}

\section{Introduction and main results} \label{sect1}

Let $W$ be the random variable of interest. Typical examples of $W$
include the
partial sum of independent random variables and functionals of
independent random variables or
dependent random variables whose joint distribution is known. Since the
exact distribution of $W$
is not available for most cases, it is natural to seek the asymptotic
distribution of $W$ with a
Berry--Esseen type error. Let $(W, W')$ be an exchangeable pair. Assume that
%
%
\begin{equation}\label{0.1}
E(W-W' | W) = g(W)+r(W),
\end{equation}
where $g(W)$ is a dominated term while $r(W)$ is a negligible term.
When $g(W) = \la W$, and $E((W'-W)^2|W)$ is concentrated around a constant,
Stein's method for normal approximation shows that
the limiting distribution of $W$ is normal
under certain regularity conditions.
We refer to \citet{Ste86},
\citet{RR97}, \citet{CS05} and references therein
for the
general theory of Stein's method.
The main aim of this paper is to find the limiting distribution of $W$
as well as
the rate of convergence for general $g$.
The key step is to identify the limiting density function. As soon as
the limiting
density function is determined, we can follow the idea of the Stein's
method of exchangeable pairs
for normal approximation. Let
%
%
\begin{equation} \label{p0}
G(t) = \int_0^t g(s) \,ds
\quad\mbox{and}\quad
p(t) = c_1 e^{-c_0 G(t)},
\end{equation}
where $c_0>0$ is a constant that will be specified later
and $c_1 = 1/\int_{-\infty}^\infty e^{-c_0 G(t)} \,dt$
is the normalizing constant. Let $Y$ be a random variable
with the probability density function $p$. Set:

\begin{enumerate}[(H1)]
\item[(H1)] $g(t)$ is nondecreasing, and $g(t)\geq0$ for $t>0$ and
$g(t) \leq0$ for $t \leq0$;

\item[(H2)] there exists $c_2 < \infty$ such that for all $x$,
%
\[
\min\bigl( 1 /c_1, 1/|c_0 g(x)|\bigr)( |x| + 3 /c_1 ) \max(1, c_0 |g'(x)|)
\leq c_2;
\]

\item[(H3)] there exists $c_3 < \infty$ such that for all $x$,
%
\[
\min\bigl( 1 /c_1, 1/|c_0 g(x)|\bigr)( |x| + 3 /c_1 ) c_0 |g'(x)|
\leq c_3.
\]

\end{enumerate}

Let $\Delta= W- W'$. Our main result shows that $W$ converges to $Y$
in distribution as long as
$c_0 E(\Delta^2 |W)$ satisfies a law of large numbers.
\begin{theorem} \label{t2.1}
Let $h$ be absolutely continuous with $\|h'\| ={\sup_x} |h'(x)|< \infty$.

\begin{longlist}[(ii)]
\item[(i)]
If \textup{(H1)} and \textup{(H2)} are satisfied, then
%
%
\begin{eqnarray}\label{t2.1a}
&&|Eh(W) - Eh(Y)|\nonumber\\
&&\qquad\leq\|h'\| \biggl\{
{ (1+c_2) \over c_1} E|1- (c_0/2)E(\Delta^2|W)| \\
&&\qquad\quad\hspace*{24pt}{}+{ 1\over2} c_0
(1+c_2) E|\Delta|^3 + { c_0 c_2 } E|r(W)|\biggr\}.
\nonumber
\end{eqnarray}
\item[(ii)] If \textup{(H1)} and \textup{(H3)} are satisfied, then
%
%
\begin{eqnarray}\label{t2.1b}
&&|Eh(W) - Eh(Y)|\nonumber\\
&&\qquad\leq\|h'\| \biggl\{
{ (1+c_3) \over c_1} E|1- (c_0/2)E(\Delta^2|W)| +{ 1\over2} c_0
(1+c_3) E|\Delta|^3 \\
&&\qquad\quad\hspace*{132pt}{} +
{ c_0 \over c_1} E\biggl( \biggl(|W| + { 3 \over c_1} \biggr)
|r(W)|\biggr) \biggr\}.
\nonumber
\end{eqnarray}
\end{longlist}
\end{theorem}

When $\Delta$ is bounded, next theorem gives a Berry--Esseen type inequality.
\begin{theorem}
\label{t2.2}
Assume that $|W- W'| \leq\delta$, where $\delta$ is a constant. If
\textup{(H1)} and \textup{(H3)} are satisfied, then
%
%
\begin{eqnarray}\label{t2.2a}\quad
&&|P(W \leq z) - P(Y \leq z)|\nonumber\\
&&\qquad\leq
3 E |1- (c_0/2)E(\Delta^2|W)| + c_1 \max(1, c_3) \delta+ 2 c_0
E|r(W)|/c_1 \\
&&\qquad\quad{} + \delta^3 c_0 \{(2+c_3/2)
E|c_0 g(W)|+c_1 c_3/2\}.\nonumber
\end{eqnarray}
\end{theorem}

We remark that $c_0$ can be chosen as follows.
In order to make the error term on the right-hand side of (\ref{t2.1a})
small, it is necessary that
$E| 1- (c_0/2)E(\Delta^2|W)| \to0$ and therefore $ E ( 1-
(c_0/2)E(\Delta^2|W))$ must be small and we
should choose $c_0$ so that $c_0 \sim2/ E(\Delta^2)$.
%
%

The paper is organized as follows. In Section \ref{sect2}, we give a concrete
application of our
general result to the magnetization of the Curie--Weiss model of
ferromagnets at the critical temperature, and show that the rate of
convergence achieves
$O(n^{-1/2})$. In Section \ref{sect3}, we focus on approximation by the
exponential distribution with an application to the
spectrum of the Bernoulli--Laplace Markov chain.
We present a general approach of Stein's method of exchangeable
pairs in Section~\ref{sect4} and postpone detailed proofs of our main
results to
Section \ref{sect5}.

\section{Curie--Weiss model}
\label{sect2}

Consider the Curie--Weiss model for $n$ spins at
temperature $T$, that is, the probability distribution on $\{-1,1\}^n$ that
puts mass
\[
Z_T^{-1}\exp\biggl(\frac{\sum_{1 \leq i<j \leq n} \sigma_i \sigma
_j}{Tn}\biggr)
\]
at $\sigma\in\{-1,1\}^n$, where $Z_T$ is the normalizing constant.
Let us fix $T=1$, which is the ``critical temperature'' for this model.
Now let
\[
W = W(\sigma) = n^{-3/4}\sum_{i=1}^n \sigma_i.
\]
This is a simple statistical mechanical model of ferromagnetic
interaction, sometimes called the Ising model on the complete graph.
For a detailed
mathematical treatment of this model, we refer to the book by
\citet{Ell85}.


Following ideas in \citet{simon73}, it was proved by Ellis
and Newman (\citeyear{ellisnewman78a}, \citeyear{ellisnewman78b}) that
as $n \rightarrow\infty$, the law of
$W$ converges to the distribution with density proportional to
$e^{-x^4/12}$. For various interesting extensions and refinements of
their results, let us refer to \citet{ellisetal80} and
\citet{papangelou89}.

Below, we present a Berry--Esseen bound for this noncentral limit
theorem obtained via Theorem \ref{t2.2}. Incidentally, Theorem \ref
{t2.2} can also be used to obtain similar error bounds for the other
limit theorems in the aforementioned papers (in particular, the
Curie--Weiss model at noncritical temperatures), but we prefer to stick
to this example only, since it is probably the most interesting and
relevant one.

Given a random element $\sigma$, construct $\sigma'$ by
choosing a coordinate $I$ at random and replacing $\sigma_I$ by
$\sigma_I'$,
where $\sigma_I'$ is generated from the conditional distribution of
$\sigma_I$
given $(\sigma_j)_{j\ne I}$. In other words,
we take one step of the Glauber dynamics. It is easy to
see that $(\sigma, \sigma')$ is an exchangeable pair.
Let $W'=W(\sigma')$. We shall show that (see Section \ref{sect5})
%
%
\begin{eqnarray}\label{3.1}
E\bigl| E(W-W'|W)- \tfrac{1}{3}n^{-3/2} W^3 \bigr| &=& O(n^{-2}),
\\
\label{3.2}
E\bigl| E\bigl((W'-W)^2|W\bigr) - 2n^{-3/2}\bigr| &=& O(n^{-2}),
\\
\label{3.3}
|W'- W| &=& O(n^{-3/4})
\end{eqnarray}
and
%
%
\begin{equation}\label{3.4}
E|W|^3 = O(1).
\end{equation}
Let us now explain roughly how we arrive at (\ref{3.1}), which is the
most important step. A simple computation shows that at any temperature,
%
\[
E(W-W'|W) = n^{-3/4}\bigl(m-\tanh(m/T)\bigr) + O(n^{-2}),
\]
where $m := n^{-1/4} W$ is the magnetization. Since $m\simeq0$ with
high probability when $T \ge1$, and $\tanh x = x - x^3/3 + O(x^5)$ for
$x\simeq0$, we see that the right-hand side in the above equation is
like $n^{-3/4}m (1-1/T)$ when $T >1$, while it is like $n^{-3/4} m^3/3$
when $T=1$. This is what distinguishes between the high temperature
regime \mbox{$T > 1$} and the critical temperature $T=1$, and this is
how we
arrive at (\ref{3.1}).

Let
\[
g(w) = \tfrac{ 1}{3} n^{-3/2} w^3,\qquad
c_0 = n^{3/2},\qquad\delta= O(n^{-3/4}).
\]
Then
\[
G_1(w) = c_0 \int_{0}^w g(t) \,dt = w^4/12.
\]
With the above information, it can be easily checked that by Theorem
\ref{t2.2}, we get the following theorem.
\begin{theorem}\label{t3.1}
Let $Y$ be a random variable with density function
\[
p(w) = c_1 e^{-w^4/12}\qquad \mbox{where }
c_1 = \frac{1}{ \int_{-\infty}^{\infty} e^{-w^4/12} \,d w}= { 2^{1/2}
\over3^{1/4} \Gamma(1/4)} .
\]
Then for all $z$,
%
%
\begin{equation}\label{t3.1a}
|P(W \leq z) - P(Y \leq z)|
\leq c n^{-1/2},
\end{equation}
where $c$ is an absolute constant.
\end{theorem}

Incidentally, after this manuscript was submitted, it was brought to
our attention that an article by \citet{eichelsbacherlowe09}
was in preparation, where the same result (Theorem \ref{t3.1}) is proved,
along the same lines as our proof. \citet{eichelsbacherlowe09} has
generalizations of Theorem \ref{t3.1} to some other mean-field models.

\section{Exponential limit with application to spectrum of the
Bernoulli--Laplace Markov chain}\label{sect3}

In this section, we focus on the exponential limit. Let
$(W, W')$ be an exchangeable pair satisfying
%
%
\begin{equation}\label{exp-01}
E(W- W'|W) = 1/c_0 + r(W),
\end{equation}
where $c_0 >0$ is a constant. Let $\Delta=W-W'$.
As a special case of Theorems \ref{t2.1} and \ref{t2.2} with a
constant function $g$, we have
\begin{theorem} \label{exp-t1}
Let $Y$ have the exponential distribution with mean $1$. Assume (\ref
{exp-01}) is satisfied.
\begin{longlist}[(ii)]
\item[(i)]
Let $h$ be absolutely continuous with $\|h'\| < \infty$. Then:
\begin{eqnarray}\label{exp-t1a}
&&|Eh(W) - Eh(Y)|\nonumber\\[-8pt]\\[-8pt]
&&\qquad\leq\|h'\| \{
E| 1 - (c_0/2) E(\Delta^2|W)|
+ c_0 E|\Delta|^3 + 3 c_0 E|W r(W)| \}. \nonumber
\end{eqnarray}

\item[(ii)] If $|\Delta| \leq\delta$ for some constant $\delta$, then
\begin{eqnarray} \label{exp-t1b}
&&|P(W \leq z) - P(Y \leq z)|\nonumber\\[-8pt]\\[-8pt]
&&\qquad\leq3 E| 1 - (c_0/2) E(\Delta^2|W)| + \delta+ 2 c_0 \delta^3 +
3 c_0 E|W r(W)|. \nonumber
\end{eqnarray}
\end{longlist}
\end{theorem}

We refer to \citet{CFR08} and \citet{PR09} for other general results
for the exponential approximation.

We now apply Theorem \ref{exp-t1} to the spectrum of the
Bernoulli--Laplace Markov chain, a simple model of diffusion,
following the work of \citet{CFR08}.
Two urns contain $n$ balls each. Initially the balls in each urn are
all of a single color, with urn 1 containing all
white balls, and urn 2 all black.
At each stage, a ball is picked at random from each urn and the two are
switched.
Let the state of the chain be the number of white balls in the urn 1.
Diaconis and Shahshahani
(\citeyear{DS87}) proved that $(n/4) \log(2n) + cn$ steps suffice for this
process to
reach equilibrium, in the sense that the total variation distance to
the stationary distribution is
at most $a e^{-dc}$ for positive universal constants $a$ and $d$. In
order to prove this, they used the fact that the
spectrum of the Markov chain consists of the numbers
%
%
\begin{equation}\label{exp-02}
\lambda_i= 1- i(2n-i+1)/n^2 \qquad\mbox{for } i=0, 1,\ldots, n,
\end{equation}
occurring with
multiplicities
\[
m_i = \pmatrix{2n\cr i} -\pmatrix{2n\cr i-1} \qquad\mbox{for }
i=0,1 ,\ldots, n.
\]
Let $I$ have distribution $P(I=i) =\pi_i$,
where
\[
\pi_i = { {2n\choose i} - {2n\choose i-1} \over
{2n\choose n}}
\]
for $0 \leq i \leq n$. Then $\lambda_I$ is a random eigenvalue chosen from
$\{\lambda_i, 0 \leq i \leq n\}$ in proportion to their
multiplicities. \citet{Hora98} proved that
$W=n \lambda_I +1$ converges in distribution to an exponential random
variable with mean $1$.

Noting that $n\lambda_i +1 = (n-i) (n+1-i)/n := \mu_i$, we can rewrite
$W= \mu_I$.
To apply Theorem \ref{exp-t1}, we construct an exchangeable pair $(W,
W')$ using a reversible Markov chain on
$\{0, 1,\ldots, n\}$ with transition probability matrix $K$ satisfying
\[
\pi(i) K(i,j) = \pi(j) K(j,i) \qquad\mbox{for all } i, j \in\{0,
1,\ldots, n\}.
\]
Given such a $K$, we
obtain the pair $(W, W')$ by letting
$W=u_I$ where $I$ is chosen from the equilibrium distribution $\pi$,
and $W'= \mu_J$ where
$J$ is determined by taking one step from state $I$ according to the
transition probability $K$.
As proved in \citet{CFR08}, we have (with
$\Delta= W-W'$)
\begin{eqnarray*}
E(\Delta|W) &=& { 1 \over2 n^2} - {n+1 \over2 n^2} I_{\{W=0\}},\qquad
%
%
E(W) = 1,\\
E(\Delta^2|W) &=& { 1 \over n^2} \quad\mbox{and}\quad
E|\Delta|^3 \leq6 n^{-5/2}.
\end{eqnarray*}
Now applying Theorem \ref{exp-t1}, we have the following theorem.
\begin{theorem} \label{exp-t2}
Let $Y$ have the exponential distribution with mean $1$ and
$h$ be absolutely continuous with $\|h'\| < \infty$. Then
%
%
\begin{equation} \label{exp-t2a}
|Eh(W) - Eh(Y)| \leq12 n^{-1/2}.
\end{equation}
\end{theorem}

As the difference between $W$ and $W'$ is large when $I$ is small,
Theorem \ref{exp-t1} does not provide a useful
Berry--Esseen type bound. However, using a completely different approach
and some heavy machinery,
\citet{CFR08} are able to show that
\[
{\sup_z }|P(W \leq z) - P(Y \leq z)| \leq C n^{-1/2},
\]
where $C$ is a universal constant.

\section{The Stein method via density approach}
\label{sect4}

Let $p$ be a strictly positive, absolutely continuous probability
density function, supported on $(a,b)$,
where $-\infty\le a < b\le\infty$.
Assume that a right limit $p(a+)$ at $a$ and a left limit $p(b-)$
exist. Let $p'$ be a version of the derivative of $p$ and assume that
\[
\int_a^b |p'(t)| \,dt < \infty.
\]
Let $Y$ be a random variable with the probability density function
$p$. In this section, we develop the Stein method via density approach.
The approach was developed in \citet{SDHR04}, but the properties
presented in Section 4.2 are new.

\subsection{The Stein identity and equation}

A key step is to have Stein's identity and Stein's equation.
Let $\mathcal{D}$ be the set of bounded, absolutely continuous
functions $f$ with
$f(b-)=f(a+)=0$. Observe that for any $f\in\mathcal{D}$
\begin{eqnarray} \label{0.3}
E\{ f'(Y) + f(Y)p'(Y)/p(Y) \}
& =& E \{(f(Y) p(Y))'/p(Y)\} \nonumber\\[-8pt]\\[-8pt]
& =& \int_{a}^b (f(y) p(y))' \,dy = 0.\nonumber
\end{eqnarray}
The Stein identity is
%
%
\begin{equation}\label{stein}
Ef'(Y) + Ef(Y) p'(Y)/p(Y) =0 \qquad\mbox{for } f \in\mathcal{D}.
\end{equation}
For any measurable function $h$ with $E|h(Y)| < \infty$,
let $f=f_h$ be the solution to Stein's
equation
%
%
\begin{equation}\label{0.4}
f'(w) + f(w)p'(w) / p(w) 
= h(w) - Eh(Y).
\end{equation}
It follows from (\ref{0.4}) that
\[
(f(w) p(w))' = \bigl(h(w) - Eh(Y)\bigr) p(w)
\]
and hence
\begin{eqnarray} \label{0.7}
f(w) &=& 1/p(w) \int_{a}^w \bigl(h(t) - Eh(Y)\bigr) p(t) \,dt
\nonumber\\[-8pt]\\[-8pt]
& =& -1/p(w) \int_{w}^{b} \bigl(h(t) - Eh(Y)\bigr) p(t) \,dt.\nonumber
\end{eqnarray}
Note that $f_h\in\mathcal{D}$.

Consider two classes of density functions. The first one is the family
of exponential distributions.
It is easy to see that if $Y$ has the exponential distribution with parameter
$\lambda$, that is, $Y$ is a random variable with density function
$p(x) = \lambda e^{-\lambda x}$
for $x>0$ and $p(x) =0$ for $ x \leq0$. Then $p'(x) /p(x) = - \lambda
$ and the Stein identity
(\ref{stein}) becomes
%
%
\begin{equation}\label{stein-exp}
Ef'(Y) - \lambda Ef(Y) =0 \qquad\mbox{for } f \in\mathcal{D}.
\end{equation}
The second is the family
\[
p(x) = { \alpha e^{-|x|^\alpha/\beta} \over2 \beta^{1/\alpha}
\Gamma(1/\alpha)}, \qquad- \infty< x < \infty,
\]
where $\alpha>0, \beta>0$. Then $p'(x)/p(x) = - { \alpha\over\beta
} |x|^{\alpha-1} \mbox{sign}(x)
$ and hence the Stein identity reduces to
\[
Ef'(Y) - { \alpha\over\beta} E|Y|^{\alpha-1} \mbox{sign}(Y) f(Y)
=0 \qquad\mbox{for } f \in\mathcal{D}.
\]

\subsection{Properties of the Stein solution}

In order to determine error bounds for the approximation to $E(h(Y))$,
we need
to understand some basic properties of the Stein solution $f_h$. In the
following, we use the notation $\|g\| := {\sup_{x\in\rr}}|g(x)|$.
\begin{lemma} \label{l3.1}
Let $h$ be a measurable function and $f_h$ be the Stein solution
and let $F(x) = \int_{a}^x p(t) \,dt$.

\begin{longlist}[(ii)]
\item[(i)]
Assume that $h$ is bounded and that there exist $d_1>0$ and $d_2>0$
%
%
\begin{equation} \label{assump-1}
\min\bigl(1-F(x), F(x)\bigr) \leq d_1 p(x)
\end{equation}
and
%
%
\begin{equation}\label{assump-2}
|p'(x)|\min\bigl(F(x), 1- F(x)\bigr)\leq d_2 p^2(x).
\end{equation}
Then
%
%
\begin{eqnarray}
\label{l3.1a}
\|f_h\| &\leq&2 d_1 \| h \|,
\\
\label{l3.1a-01}
\| f_h p'/p\| &\leq&2 d_2 \|h\|
\end{eqnarray}
and
%
%
\begin{equation}\label{l3.1b}
\|f'_h \|\leq(2 + 2 d_2) \|h\|.
\end{equation}
\item[(ii)] Assume that $h$ is absolutely continuous with bounded $h'$.
In addition to (\ref{assump-1}), (\ref{assump-2}), assume that there exist
$d_3$ and $d_4$ such that
\begin{eqnarray}\label{assump-3}\quad
&&\min\bigl( E|Y|I_{\{Y \leq x\}} + E|Y| F(x),
E|Y|I_{\{Y > x\}} + E|Y|\bigl(1- F(x)\bigr) \bigr)
|(p'/p)'|\nonumber\\[-8pt]\\[-8pt]
&&\qquad\leq d_3 p(x)
\nonumber
\end{eqnarray}
and
\begin{eqnarray}\label{assump-4}
&&\min\bigl( E|Y|I_{\{Y \leq x\}} + E|Y| F(x),
E|Y|I_{\{Y > x\}} + E|Y|\bigl(1- F(x)\bigr) \bigr)\nonumber\\[-8pt]\\[-8pt]
&&\qquad\leq d_4 p(x).
\nonumber
\end{eqnarray}
Then if $h$ is absolutely continuous with bounded derivative $h'$,
%
%
\begin{eqnarray}
\label{l3.1c}
\|f_h''\| &\leq&(1+ d_2)(1+d_3) \|h'\|,
\\
\label{l3.1d}
\|f_h \| &\leq& d_4 \|h'\|
\end{eqnarray}
and
%
%
\begin{equation}\label{l3.1e}
\|f_h' \| \leq(1+d_3) d_1 \|h'\|.
\end{equation}
\end{longlist}
\end{lemma}
\begin{pf} (i) Let $Y^*$ be an independent copy of $Y$. Then we can
rewrite $f_h$ in (\ref{0.7}) as
\begin{eqnarray}\label{l3.1-0}
f(w) &=& \bigl(1/p(w)\bigr) E\bigl(h(Y) - h(Y^*)\bigr)I_{\{Y \leq w\}},
\nonumber\\[-8pt]\\[-8pt]
&=&-\bigl(1/p(w)\bigr) E\bigl(h(Y) - h(Y^*)\bigr)I_{\{Y > w\}}
,\nonumber
\end{eqnarray}
which yields
%
%
\begin{equation} \label{l3.1-1}
|f(w)| \leq2 \|h\| \min\bigl(F(w), 1- F(w)\bigr) /p(w).
\end{equation}
Inequality (\ref{l3.1a}) now follows from (\ref{assump-1}) and (\ref
{l3.1-1}).
Inequalities (\ref{l3.1-1}) and (\ref{assump-2}) imply
$|f_h p'/p| \leq2 d_2 \|h\|$, that is (\ref{l3.1a-01}), and now
(\ref{l3.1b}) follows from (\ref{0.4}).

(ii) Let $g_1(x) = p'(x)/p(x)$. Recall by (\ref{0.4})
%
%
\begin{equation}\label{l3.1-2}
f'' = h' - f' g_1 - f g_1'.
\end{equation}
To prove (\ref{l3.1c}), it suffices to show that
%
%
\begin{equation}\label{l3.1-3}
\|f g_1'\| \leq d_3 \| h'\|
\end{equation}
and
%
%
\begin{equation}\label{l3.1-4}
\|f' g_1\|
\leq(1+d_3) d_2 \| h'\|.
\end{equation}
By (\ref{l3.1-0}) again, we have
%
%
\begin{eqnarray} \label{l3.1-5}\qquad
|f(w) p(w) |
& \leq& \|h'\| \min\bigl( E(|Y|+ |Y^*|)I_{\{Y \leq w\}},
E(|Y|+ |Y^*|)I_{\{Y > w\}}\bigr)\nonumber\\
& = & \|h'\| \min\bigl( E|Y|I_{\{Y \leq w\}} + E|Y| F(w),
E|Y|I_{\{Y > w\}}\\
&&\hspace*{127.2pt}{} + E|Y|\bigl(1- F(w)\bigr)\bigr).\nonumber
\end{eqnarray}
This proves (\ref{l3.1-3}) by assumption (\ref{assump-2}). This also proves
(\ref{l3.1d})
by (\ref{assump-4}).

It follows from (\ref{l3.1-2}) that
\[
(h'- f g_1')p = p( f'' + f' g_1) = f'' p + f' p' = (f'p)'.
\]
Thus
\[
f'(w) p(w) = \int_a^w (h'- f g_1')p\, dx
= - \int_w^b (h'- fg_1') p \,dx
\]
and hence
\[
|f'(w) p(w)| \leq\|h'\| ( 1+ d_3) \min\bigl( F(w), 1-F(w)\bigr),
\]
which gives (\ref{l3.1-4}) as well as (\ref{l3.1e}) by (\ref
{assump-4}) and
(\ref{assump-1}), respectively.
\end{pf}

The next lemma shows that (\ref{assump-1})--(\ref{assump-4}) are satisfied
for $p$ defined in (\ref{p0}).
\begin{lemma}
\label{l3.2}
Let $p$ be defined as in (\ref{p0}). Assume that \textup{(H1)} and
\textup{(H2)} are
satisfied. Then
(\ref{assump-1})--(\ref{assump-4}) hold with
$d_1= 1/c_1$, $d_2=1$, $d_3= c_2$ and $d_4=c_2$.
\end{lemma}
\begin{pf}
Let $g_2(t)=c_0 g(t)$, $G_1(t)= c_0 G(t)$ and
$F(t) = P(Y \leq t)$ be the distribution function of $Y$.
We first show that (\ref{assump-1}) is satisfied with $d_1=1/c_1$. It
suffices to show that
%
%
\begin{equation}\label{l3.2-1}
F(t) \leq F(0) p(t)/c_1 \qquad\mbox{for } t \leq0
\end{equation}
and
%
%
\begin{equation}\label{l3.2-2}
1-F(t) \leq\bigl(\bigl(1-F(0)\bigr)/c_1\bigr) p(t)
\qquad\mbox{for } t \geq0.
\end{equation}
Let $H(t)=F(t) - (F(0)/c_1)p(t)$ for $ t\leq0$.
Noting that
\begin{eqnarray*}
H'(t) &=& p(t) - \bigl(F(0)/c_1\bigr)p'(t)\\
& = & p(t) + \bigl(F(0)/c_1\bigr) g_2(t) p(t) \\
& =& p(t) \bigl( 1+ g_2 (t) F(0) /c_1\bigr).
\end{eqnarray*}
Since $g_2(t)$ is nondecreasing, if $H'(0)>0$, then there is at most
one $t_0$ such that
$H'(t_0)=0$; if $H'(0) \leq0$, then $H'(t) \leq0$ for $ t < 0$. Hence,
$H$ achieves maximum either at $t=0$ or $t=-\infty$. Notice that
$H(0)= H(-\infty) =0$, $H(t) \leq0$ for all $ t < 0$. This proves
(\ref{l3.2-1}). Similarly,
(\ref{l3.2-2}) holds.

Next, we prove (\ref{assump-2}). Noting that
$p'= -p g_2$, we have for $t<0$
\begin{eqnarray}\label{l3.2-3}
F(t) & =& \int_{-\infty}^t p(s) \,ds\nonumber\\
& \leq& \int_{-\infty}^t { g_2(s) p(s) \over g_2(t)} \,ds
\nonumber\\[-8pt]\\[-8pt]
& =&
\int_{-\infty}^t {- p'(s) \over g_2(t)} \,ds \nonumber\\
& =& { p(t) \over- g_2(t)} = p(t) /|g_2(t)| . \nonumber
\end{eqnarray}
Similarly, we have
%
%
\begin{equation}\label{l3.2-4}
1-F(t) \leq p(t)/g_2(t)
\qquad\mbox{for } t \geq0.
\end{equation}
Hence, (\ref{assump-2}) is satisfied with $d_2=1$.

Note that (\ref{assump-1}) and (\ref{assump-2}) imply that
%
%
\begin{equation}\label{l3.2-5}
1-F(x)
\leq p(x) \min\bigl( 1 /c_1, 1/|g_2(x)|\bigr)
\qquad\mbox{for } x\geq0
\end{equation}
and
%
%
\begin{equation}\label{l3.2-6}
F(x)\leq p(x) \min\bigl( 1 /c_1, 1/|g_2(x)|\bigr)
\qquad\mbox{for } x\leq0.
\end{equation}

To verify (\ref{assump-3}), with $x\geq0$ write
\begin{eqnarray}\label{l3.2-7}
E|Y|I_{\{Y > x\}}
& =& x P(Y > x) + \int_x^\infty P(Y\geq t) \,dt\nonumber\\
& \leq& x p(x) \min\bigl( 1 /c_1, 1/|g_2(x)|\bigr)\nonumber\\
&&{} + \int_x^\infty p(t) \min\bigl( 1 /c_1, 1/|g_2(t)|\bigr) \,dt
\nonumber
\\
& \leq& x p(x) \min\bigl( 1 /c_1, 1/|g_2(x)|\bigr)\nonumber\\[-8pt]\\[-8pt]
&&{} + \min\bigl( 1 /c_1, 1/|g_2(x)|\bigr) \int_x^\infty p(t) \,dt
\nonumber
\\
& \leq& \min\bigl( 1 /c_1, 1/|g_2(x)|\bigr)\bigl\{ x p(x) + \bigl(1-
F(x)\bigr)\bigr\}
\nonumber\\
& \leq& \min\bigl( 1 /c_1, 1/|g_2(x)|\bigr)\{ x p(x) + p(x)/c_1
\} \nonumber\\
&\leq& p(x) \min\bigl( 1 /c_1, 1/|g_2(x)|\bigr)\{ x + 1 /c_1 \}.
\nonumber
\end{eqnarray}
Similarly, for $x< 0$,
%
%
\begin{equation}\label{l3.2-8}
E|Y|I_{\{Y < x\}}
\leq p(x) \min\bigl( 1 /c_1, 1/|g_2(x)|\bigr)\{ |x| + 1 /c_1 \}.
\end{equation}
Equations (\ref{l3.2-7}) and (\ref{l3.2-8}) with $x=0$ also give
$E|Y| \leq2/c_1$. Hence, recalling (\ref{l3.2-5})
\begin{eqnarray} \label{l3.2-9}
&&E|Y|I_{\{Y > x\}} + E|Y| \bigl(1-F(x)\bigr)
\nonumber\\[-8pt]\\[-8pt]
&&\qquad\leq p(x) \min\bigl( 1 /c_1, 1/|g_2(x)|\bigr)\{ x + 3 /c_1 \}
\qquad\mbox{for } x>0
\nonumber
\end{eqnarray}
and
\begin{eqnarray} \label{l3.2-10}
&&E|Y|I_{\{Y < x\}} + E|Y| F(x)
\nonumber\\[-8pt]\\[-8pt]
&&\qquad\leq p(x) \min\bigl( 1 /c_1, 1/|g_2(x)|\bigr)\{ |x| + 3 /c_1 \}
\qquad\mbox{for } x\leq0.
\nonumber
\end{eqnarray}
Thus, (\ref{assump-3}) holds with $d_3= c_2$ by (H2).

Equations (\ref{l3.2-9}) and (\ref{l3.2-10}) also show that (\ref{assump-4})
is satisfied with $d_4=c_2$.

This completes the proof of Lemma \ref{l3.2}.
\end{pf}

From the proof of Lemma \ref{l3.2}, one can see the following remark
is true.
\begin{remark} \label{r3.1}
Assume that (H1) and (H3) are satisfied. Then (\ref
{assump-1})--(\ref{assump-3}) hold with
$d_1=1/c_1$, $d_2=1$ and $d_3=c_3$, and hence (\ref{l3.1c}) and (\ref{l3.1e}).
\end{remark}

\section{Proof of main results}
\label{sect5}

In this section, we prove the general error bounds (Theorems \ref
{t2.1} and \ref{t2.2}), the result for the Curie--Weiss model
(Theorem \ref{t3.1}), and Theorem \ref{exp-t1}.

\subsection{\texorpdfstring{Proof of Theorem \protect\ref{t2.1}}{Proof of Theorem 1.1}}
Let $f= f_h$ be the solution to Stein's equation (\ref{0.4}). Then
\begin{eqnarray} \label{0.20}
Eh(W) - Eh(Y) &= & Ef'(W)+ Ef(W) p'(W)/ p(W)\nonumber\\[-8pt]\\[-8pt]
&= & Ef'(W) - c_0 Ef(W)g(W).
\nonumber
\end{eqnarray}
Recall $\Delta= W-W'$ and observe that
for any absolutely continuous function $f$
\begin{eqnarray} \label{0.2}
0 & =& E(W-W')\bigl( f(W') + f(W)\bigr) \nonumber\\
& = & 2 Ef(W) (W-W') + E(W-W')\bigl( f(W') -
f(W)\bigr)\nonumber\\[-8pt]\\[-8pt]
& =& 2 E\bigl\{f(W) E\bigl((W-W')|W\bigr)\bigr\}
- E(W-W')\int_{-\Delta}^0 f'(W+t) \,dt \nonumber\\
& =& 2 Ef(W) g(W) + 2 Ef(W) r(W) - E\int_{-\infty}^\infty f'(W+t)
\hatK(t) \,dt,\nonumber
\end{eqnarray}
where
\[
\hatK(t) = E\bigl\{\Delta( I\{- \Delta\leq t \leq0\} -
I\{0 < t \leq-\Delta\})|W\bigr\}.
\]
Substituting (\ref{0.2}) into (\ref{0.20}) gives
%
%
\begin{eqnarray} \label{0.6}
&&Ef'(W) - c_0 Ef(W)g(W)\nonumber\\
&&\qquad= Ef'(W) - (c_0/2) \biggl\{ E\int_{-\infty}^\infty f'(W+t)
\hatK
(t) \,dt - 2Ef(W) r(W)\biggr\}
\nonumber\\
&&\qquad= E\bigl\{f'(W)\bigl( 1- (c_0/2) E(\Delta^2|W)\bigr)\bigr\}\\
&&\qquad\quad{} + (c_0/2) E\int_{-\infty}^\infty\bigl(f'(W) -
f'(W+t)\bigr)\hatK(t) \,dt
\nonumber
\\
&&\qquad\quad{} + c_0 Ef(W)r(W).\nonumber
\end{eqnarray}
When (H1) and (H2) are satisfied, by Lemmas \ref{l3.1} and \ref{l3.2}
%
%
\begin{equation} \label{t2.1-0}
\|f_h\| \leq c_2 \|h'\|,\qquad
\|f_h'\| \leq(1+c_2) \|h'\|/c_1,\qquad
\|f_h''\| \leq2 (1+c_2) \|h'\|\hspace*{-28pt}
\end{equation}
and hence
\begin{eqnarray*}
&&|Ef_h'(W) - c_0 Ef_h(W)g(W)|\\
&&\qquad\leq{ (1+c_2) \|h'\| \over c_1} E\bigl|\bigl( 1- (c_0/2)
E(\Delta^2|W)\bigr)\bigr|
\\
&&\qquad\quad{} + (1+c_2) \|h'\| c_0 E|\Delta|^3 /2 + c_0 c_2 \|h'\| E|r(W)|.
\end{eqnarray*}
This proves (\ref{t2.1a}).

Under (H1) and (H3), by Remark \ref{r3.1}
%
%
\begin{equation} \label{t2.1-001}
\|f_h'\| \leq(1+c_3) \|h'\|/c_1,\qquad
\|f_h''\| \leq2 (1+c_3) \|h'\|.
\end{equation}
From (\ref{l3.1-0}), (\ref{l3.2-9}) and (\ref{l3.2-10}) it follows that
%
%
\begin{eqnarray} \label{t2.1-002}\quad
|f(w)| &\leq& \bigl(1/p(w)\bigr) \|h' \| \min\bigl( E|Y-Y*| I_{\{Y \leq
w\}},
E|Y-Y*| I_{\{Y \geq w\}}\bigr)\nonumber\\
& \leq& \| h'\| \min\bigl( 1/c_1, 1/|g_2(w)|\bigr) (|w| + 3/c_1)\\
&\leq&\|h'\| (|w| + 3 /c_1) /c_1.\nonumber
\end{eqnarray}
This proves (\ref{t2.1b}) by (\ref{0.6}), (\ref{t2.1-001}) and (\ref
{t2.1-002}).

\subsection{\texorpdfstring{Proof of Theorem \protect\ref{t2.2}}{Proof of Theorem 1.2}}
Since (\ref{t2.2a}) is trivial when
$c_1 c_3 \delta>1$, we assume
%
%
\begin{equation}\label{t2.2-00}
c_1 c_3 \delta\leq1.
\end{equation}
Let $F$ be the distribution function of $Y$ and let $f=f_z$ be the
solution to the equation
%
%
\begin{equation}\label{t2.2-1}
f'(w) - c_0 f(w) g(w) = I(w \leq z) - F(z).
\end{equation}
By (\ref{0.2}),
\begin{eqnarray*}
&&2 Ef(W)g(W) + 2Ef(W) r(W)\\
&&\qquad= E
\int_{-\infty}^\infty f'(W+t) \hatK(t) \,dt
\\
&&\qquad= E\int_{-\delta}^\delta\{c_0 f(W+t) g(W+t) + I(W+t \leq z) -
F(z)\}\hatK(t) \,dt \\
&&\qquad\geq E \int_{-\delta}^\delta c_0 f(W+t) g(W+t)\hatK(t) \,dt +
EI(W\leq z-\delta) \Delta^2 - F(z) E\Delta^2
\end{eqnarray*}
and hence
\begin{eqnarray} \label{t2.2-30}
&&EI(W\leq z-\delta) \Delta^2 - F(z) E\Delta^2\nonumber\\
&&\qquad\leq2Ef(W)g(W) +
2Ef(W) r(W) \nonumber\\
&&\qquad\quad{} - c_0 E\int_{-\delta}^\delta f(W+t) g(W+t) \hatK(t)
\,dt\nonumber\\[-8pt]\\[-8pt]
&&\qquad= 2Ef(W)g(W)\bigl( 1- (c_0/2)E(\Delta^2|W)\bigr) + 2 Ef(W)
r(W) \nonumber\\
&&\qquad\quad{} + c_0 E\int_{-\delta}^\delta\{f(W) g(W) - f(W+t)
g(W+t)\} \hatK
(t) \,dt\nonumber\\
&&\hspace*{-3.4pt}\qquad:= J_1 + J_2 + J_3.\nonumber
\end{eqnarray}
From Lemmas \ref{l3.1} and \ref{l3.2} again, we obtain
%
%
\begin{equation}\label{t2.2-0}
\|f_z\| \leq2/c_1,\qquad\|f_z g\| \leq2 /c_0 \quad\mbox{and}\quad
\|f_z'\| \leq4.
\end{equation}
Therefore,
%
%
\begin{equation} \label{t2.2-32}
|J_1| \leq(4/c_0) E| 1- (c_0/2)E(\Delta^2|W)|.
\end{equation}
and
%
%
\begin{equation} \label{t2.2-31}
|J_2|\leq(4/c_1) E|r(W)|.
\end{equation}
To bound $J_3$, we first show that
%
%
\begin{equation}\label{t2.2-01}
{\sup_{|t|\leq\delta}} |g(w+t) - g(w)|
\leq{ c_1 c_3 \delta\over2 c_0} \bigl( c_1 + c_0 |g(w)|\bigr).
\end{equation}
From (H2), it follows that
%
%
\begin{eqnarray}\label{t2.2-02}
|g'(x)| &\leq& { c_1 c_3 \over3 c_0 \min( 1/c_1, 1/|c_0
g(x)|)}\nonumber\\
& =&{ c_1 c_3 \over3 c_0} \max( c_1, |c_0 g(x)|) \\
& \leq& { c_1 c_3 \over3 c_0} \bigl( c_1 + |c_0 g(x)|\bigr). \nonumber
\end{eqnarray}
Thus, by the mean value theorem,
\begin{eqnarray*}
&&{\sup_{|t|\leq\delta}} |g(w+t) - g(w)|\\
&&\qquad \leq {\delta\sup_{|t|\leq\delta} }|g'(w+t)| \\
&&\qquad \leq { c_1 c_3 \delta\over3 c_0} \Bigl( c_1 + {c_0 \sup_{|t| \leq
\delta}}
|g(w+t)|\Bigr)\\
&&\qquad \leq { c_1 c_3 \delta\over3 c_0} \Bigl( c_1 + c_0 |g(w)| + {c_0 \sup
_{|t| \leq\delta}}
|g(w+t)- g(w)|\Bigr) \\
&&\qquad = { c_1 c_3 \delta\over3 c_0} \bigl( c_1 + c_0 |g(w)|\bigr)
+ {{ c_1 c_3 \delta\over3} \sup_{|t| \leq\delta}}
|g(w+t)- g(w)| \\
&&\qquad \leq { c_1 c_3 \delta\over3 c_0} \bigl( c_1 + c_0 |g(w)|\bigr)
+ {{1 \over3} \sup_{|t| \leq\delta}} |g(w+t)- g(w)|
\end{eqnarray*}
by (\ref{t2.2-00}). This proves (\ref{t2.2-01}).

Now by (\ref{t2.2-0}) and (\ref{t2.2-01}), when $|t|\leq\delta$
\begin{eqnarray*}
&&|f(w) g(w) - f(w+t) g(w+t)|\\
&&\qquad \leq |g(w)||f(w+t) - f(w)| + |f(w+t)||g(w+t) - g(w)| \\
&&\qquad \leq 4|g(w)| |t| + { 2\over c_1} { c_1 c_3 \delta\over2 c_0}
\bigl( c_1 + c_0 |g(w)|\bigr)\\
&&\qquad \leq (4+c_3) \delta|g(w)| + \delta c_1 c_3/c_0.
\end{eqnarray*}
Therefore,
%
%
\begin{eqnarray}\label{t2.2-33}
|J_3| &\leq& c_0 (4+ c_3)\delta E|g(W)| \Delta^2 + \delta c_1 c_3
E\Delta^2
\nonumber\\[-8pt]\\[-8pt]
&\leq& (4+ c_3)\delta^3 E|c_0 g(W)| + c_1 c_3 \delta^3.\nonumber
\end{eqnarray}
Combining (\ref{t2.2-30}), (\ref{t2.2-31}), (\ref{t2.2-32}) and
(\ref{t2.2-33})
shows that
%
%
\begin{eqnarray}\label{t2.2-34}
&&EI(W\leq z-\delta) \Delta^2 - F(z) E\Delta^2\nonumber\\
&&\qquad \leq (4/c_0) E| 1- (c_0/2)E(\Delta^2|W)| + (4/c_1) E|r(W)|
\\
&&\qquad\quad{} + (4+ c_3)\delta^3 E|c_0 g(W)| + c_1 c_3 \delta^3
.\nonumber
\end{eqnarray}
On the other hand, using $F'(z) = p(z) \leq c_1$, we have
%
%
\begin{eqnarray}\label{t2.2-34-0}
&&EI(W\leq z-\delta) \Delta^2 - F(z) E\Delta^2\nonumber\\
&&\qquad ={ 2 \over c_0} \bigl( EI(W \leq z-\delta) - F(z-\delta)\bigr)
\nonumber\\
&&\qquad\quad{} - { 2 \over c_0} E \biggl\{ \bigl(I( W \leq z-\delta) - F(z)\bigr)\biggl(1- { c_0
\over2} E(\Delta^2 |W)\biggr)\biggr\} \nonumber\\[-8pt]\\[-8pt]
&&\qquad\quad{} + { 2 \over c_0} \bigl(F(z-\delta) - F(z)\bigr) \nonumber\\
&&\qquad \geq { 2 \over c_0} \bigl( P( W \leq z- \delta) - F(z-\delta)
\bigr) \nonumber\\
&&\qquad\quad{} - { 2 \over c_0} E \biggl| 1- { c_0\over2} E(\Delta^2|W)\biggr| - { 2 c_1
\delta\over c_0},\nonumber
\end{eqnarray}
which together with (\ref{t2.2-34}) yields
%
%
\begin{eqnarray} \label{t2.2-35}
&&P(W \leq z-\delta) - F(z-\delta)\\
&&\qquad \leq E |1- (c_0/2)E(\Delta^2|W)| + c_1 \delta\nonumber\\
&&\qquad\quad{} + { c_0 \over2} \bigl( (4/c_0) E| 1- (c_0/2)E(\Delta^2|W)| +
(4/c_1) E|r(W)| \nonumber\\
&&\qquad\quad\hspace*{93.7pt}{} + (4+ c_3)\delta^3 E|c_0 g(W)| + c_1 c_3 \delta^3 \bigr) \nonumber
\\
&&\qquad = 3 E |1- (c_0/2)E(\Delta^2|W)| + c_1 \delta+ 2 c_0
E|r(W)|/c_1\nonumber\\
&&\qquad\quad{} + \delta^3 c_0 \{ ( 2 + c_3 /2) E|c_0 g(W)| + c_1 c_3 /2\}.
\end{eqnarray}
Similarly, we have
%
%
\begin{eqnarray} \label{t2.2-36}
&&F(z+\delta) - P(W \leq z+\delta)\\
&&\qquad \leq 3 E |1- (c_0/2)E(\Delta^2|W)| + c_1 \delta+ 2 c_0
E|r(W)|/c_1\nonumber\\
&&\qquad\quad{} + \delta^3 c_0 \{ ( 2 + c_3 /2) E|c_0 g(W)| + c_1 c_3 /2\}.
\end{eqnarray}
This completes the proof of (\ref{t2.2a}).

\subsection{\texorpdfstring{Proof of Theorem \protect\ref{t3.1}}{Proof of Theorem 2.1}}

By (\ref{3.1})--(\ref{3.4})
\begin{eqnarray*}
E|r(W)|&=&O(n^{-2}),
\\
E\bigl| 1 - (c_0/2) E\bigl((W-W')^2|W\bigr)\bigr|
&=&O(n^{-1/2}),
\\
E|W|^3 &=&O(1).
\end{eqnarray*}
Applying Theorem \ref{t2.2} gives Theorem \ref{t3.1}.

We now show that (\ref{3.1})--(\ref{3.4}) hold.
\begin{lemma}
\label{l5.1}
With $W,W'$ as in Section \ref{sect2}, we have
%
%
\begin{eqnarray}\label{eq1}
E\biggl|E(W-W'|W) - \frac{n^{-3/2}}{3}W^3\biggr| &\le&15 n^{-2},
\\
\label{eq2}
E\bigl|E\bigl((W- W')^2|W\bigr) - 2n^{-3/2}\bigr| &\le&15 n^{-2}
\end{eqnarray}
and
%
%
\begin{equation} \label{eq3}
E|W|^3 \leq15.
\end{equation}
Also, obviously, $|W-W'|\le2n^{-3/4}$.
\end{lemma}
\begin{pf}
Let $m = n^{-1}\sum_{i=1}^n \sigma_i = n^{-1/4} W$, and for each $i$, let
\[
m_i = n^{-1} \sum_{j\ne i} \sigma_j.
\]
It is easy to see that for $\tau\in\{-1, 1\}$
%
%
\begin{equation}\label{0}
P(\sigma_i' = \tau|\sigma) = \frac{e^{m_i \tau}}{e^{m_i} + e^{-m_i}},
\end{equation}
and so
\[
E(\sigma_i' | \sigma) = { e^{m_i} \over e^{m_i} + e^{-m_i}}
- { e^{-m_i} \over e^{m_i} + e^{-m_i}} = \tanh m_i.
\]
%
Hence,
%
%
\begin{eqnarray}\label{exp0}
E(W - W'|\sigma) &= &\frac{1}{n}\sum_{i=1}^nn^{-3/4}\bigl(\sigma_i -
E(\sigma_i'|\sigma) \bigr)
\nonumber\\[-8pt]\\[-8pt]
&=& n^{-3/4}m - n^{-7/4} \sum_{i=1}^n \tanh m_i .\nonumber
\end{eqnarray}
Now it is easy to verify that the function
\[
\frac{d^2}{dx^2}\tanh x = \frac{-2\sinh x}{\cosh^3 x} = -2(\tanh x)
( 1- \tanh^2 x)
\]
has exactly two extrema
$\pm x^*$ on the real line, where $x^*$
solves the equation $\tanh^2 x^* = \frac{1}{3}$.
It follows that the maximum magnitude of this function is $4/3^{3/2}$.
Thus, for all $x,y\in\rr$,
\[
|{\tanh x} - \tanh y - (x-y)(\cosh y)^{-2}|\le\frac{2(x-y)^2}{3^{3/2}}.
\]
It follows that
\[
\Biggl|\sum_{i=1}^n \tanh m_i - n\tanh m + n^{-1}(\cosh m)^{-2}\sum
_{i=1}^n \sigma_i\Biggr|\le\frac{2n^{-1}}{3^{3/2}},
\]
and therefore
\[
\Biggl|\sum_{i=1}^n \tanh m_i - n\tanh m\Biggr|\le|m| + \frac
{2n^{-1}}{3^{3/2}}.
\]
Using this in (\ref{exp0}) and the relation $m = n^{-1/4} W$, we get
%
%
\begin{equation}\label{exp1}\quad
|E(W - W' |\sigma) + n^{-3/4} \tanh m - n^{-3/4}m| \le n^{-2}|W|
+ \frac{2n^{-11/4}}{3^{3/2}}.
\end{equation}
Now consider the function $f(x) = \tanh x - x + \frac{x^3}{3}$. Note
that $f'(x) = (\cosh x)^{-2} - 1 + x^2 \ge0$ for all $x$, and hence
$f$ is an increasing function. Also $f(0) = 0$. Therefore, $f(x) \ge0$
for all $x\ge0$. Now, it can be easily verified that the first four
derivatives of $f$ vanish at zero, and for all $x\ge0$,
\[
\frac{d^5f}{dx^5} = \frac{16}{\cosh^2 x} - 120 \frac{\sinh^2
x}{\cosh^4 x} + 120 \frac{\sinh^4 x}{\cosh^6 x} \le\frac
{16}{\cosh^2 x} \le16.
\]
Thus, for all $x\ge0$,
\[
0\le f(x)\le\frac{16}{5!} x^5 = \frac{2x^5}{15}.
\]
Since $f$ is an odd function, we get that for all $x$,
\[
\biggl|\tanh x - x + \frac{1}{3}x^3\biggr|\le\frac{2|x|^5}{15}.
\]
Using this information in (\ref{exp1}), we get
\[
\biggl|E(W-W'|\sigma) - \frac{n^{-3/4}}{3} m^3\biggr| \le\frac
{2n^{-3/4} |m|^5}{15} + n^{-2}|W| + \frac{2n^{-11/4}}{3^{3/2}}.
\]
Using the relation $m = n^{-1/4}W$, we get
%
%
\begin{equation}\label{exp2}\quad
\biggl|E(W-W'|\sigma) - \frac{n^{-3/2}}{3} W^3\biggr| \le\frac
{2n^{-2} |W|^5}{15}+ n^{-2}|W| + \frac{2n^{-11/4}}{3^{3/2}}.
\end{equation}
This implies, in particular, that
%
%
\begin{eqnarray}\label{moment1}
&&\biggl|E\bigl((W-W') W^3\bigr) - \frac{n^{-3/2}}{3}E(W^6)\biggr|
\nonumber\\[-8pt]\\[-8pt]
&&\qquad\le \frac{2n^{-2} E(W^8)}{15}+ n^{-2}E(W^4) + \frac
{2n^{-11/4}E|W|^3}{3^{3/2}}. \nonumber
\end{eqnarray}
Thus,
%
%
\begin{eqnarray} \label{moment2}
E(W^6) &\le& 3n^{3/2}\bigl|E\bigl((W'-W)W^3\bigr)\bigr| + \frac{2n^{-1/2}
E(W^8)}{5}\nonumber\\[-8pt]\\[-8pt]
& &{} + 3n^{-1/2}E(W^4) + \frac{2n^{-5/4}E|W|^3}{3^{1/2}}.\nonumber
\end{eqnarray}
Using the crude bound $|W|\le n^{1/4}$, we get
%
%
\begin{eqnarray}\label{moment3}
&&\frac{2n^{-1/2} E(W^8)}{5} + 3n^{-1/2}E(W^4) + \frac
{2n^{-5/4}E|W|^3}{3^{1/2}}
\nonumber\\[-8pt]\\[-8pt]
&&\qquad\le \frac{2 E(W^6)}{5}
+ 3E(W^2) + \frac{2n^{-1}E(W^2)}{3^{1/2}}. \nonumber
\end{eqnarray}
Next, note that by the exchangeability of $(W,W')$,
\begin{eqnarray*}
E\bigl((W'-W)W^3\bigr) &=& \tfrac{1}{2}E\bigl((W'-W)(W^3 - W'^3)\bigr)\\
&= & -\tfrac{1}{2}E\bigl((W'-W)^2 (W^2 + WW' + W'^2)\bigr).
\end{eqnarray*}
Since $|W- W'|\le2n^{-3/4}$, this gives
%
%
\begin{equation}\label{moment4}
\bigl|E\bigl((W'-W)W^3\bigr)\bigr|\le6n^{-3/2} E(W^2).
\end{equation}
Combining (\ref{moment2}), (\ref{moment3}) and (\ref{moment4}), we get
\[
E(W^6) \le\biggl(21 + \frac{2n^{-1}}{3^{1/2}}\biggr) E(W^2) + \frac
{2E(W^6)}{5},
\]
and therefore,
\[
E(W^6) \le\frac{5}{3} \biggl(21 + \frac{2n^{-1}}{3^{1/2}}\biggr)
E(W^2)\le36.9245 E(W^2).
\]
Since $E(W^2) \le(E(W^6))^{1/3}$, this gives
%
%
\begin{equation}\label{moment5}
E(W^6)\le(36.9245)^{3/2}\le224.4
\end{equation}
and hence (\ref{eq3}) holds.

Combined with (\ref{exp2}), this gives
%
%
\begin{eqnarray} \label{exp3}
&&E\biggl|E(W-W'|W) - \frac{n^{-3/2}}{3}W^3\biggr|\nonumber\\[-8pt]\\[-8pt]
&&\qquad\le n^{-2}\biggl(\frac{2(224.4)^{5/6}}{15} + (224.4)^{1/6}\biggr)
+ \frac{2n^{-11/4}}{3^{3/2}} \le15 n^{-2}.
\nonumber
\end{eqnarray}
By (\ref{0}), we have
\begin{eqnarray*}
E\bigl((W-W')^2|\sigma\bigr) &= &\frac{1}{n}\sum_{i=1}^n 4n^{-3/2}\frac
{e^{-m_i\sigma_i}}{e^{m_i\sigma_i} + e^{-m_i\sigma_i}}\\
&=& 2n^{-5/2}\sum_{i=1}^n \bigl(1-\tanh(m_i\sigma_i)\bigr)\\
&=& 2n^{-3/2} - 2n^{-5/2}\sum_{i=1}^n \sigma_i \tanh m_i.
\end{eqnarray*}
Using $|{\tanh m_i} - \tanh m |\le|m_i - m|\le n^{-1}$, we get
\begin{eqnarray*}
\bigl|E\bigl((W-W')^2 |\sigma\bigr) - 2n^{-3/2}\bigr|
&\le& 2n^{-5/2} + 2n^{-3/2}m\tanh m \\
&\le& 2n^{-5/2} + 2n^{-3/2} m^2\\
&= & 2n^{-5/2} + 2n^{-2} W^2.
\end{eqnarray*}
Using (\ref{moment5}), we get
\[
E\bigl|E\bigl((W- W')^2|W\bigr) - 2n^{-3/2}\bigr|
\le2n^{-5/2} + 2n^{-2}(224.4)^{1/3} \le15 n^{-2}.
\]
This completes the proof of the lemma.
\end{pf}

\subsection{\texorpdfstring{Proof of Theorem \protect\ref{exp-t1}}{Proof of Theorem 3.1}}
With $p(w) = e^{-w}I_{\{w >0\}}$, for given $h$, let $f_h$ be the Stein
solution given in (\ref{0.7})
\[
f_h(w) = e^{w} \int_0^w \bigl(h(t) - Eh(Y)\bigr) e^{-t} \,dt
= - e^{w} \int_{w}^\infty\bigl(h(t) - Eh(Y)\bigr)e^{-t} \,dt
\]
for $w \geq0$. Following the proof of Theorems \ref{t2.1} and \ref
{t2.2}, it suffices to show that
%
%
\begin{equation}\label{exp-t1-1}
|f_h(w)| \leq3 \min(\|h\|, \|h'\|) w \qquad\mbox{for } w \geq0.
\end{equation}
By (\ref{l3.1-1}),
\[
|f_h(w)| \leq2 \|h\| \min( 1- e^{-w}, e^{-w}) e^{w}
= 2 \|h\| \min( 1, e^w -1) \leq3 w \|h\|
\]
and by (\ref{l3.1-5})
\begin{eqnarray*}
|f_h(w)| & \leq& \|h'\| e^w
\min\bigl( - w e^{-w} + 2 (1- e^{-w}), (w+1)e^{-w}\bigr) \\
& \leq& \|h'\| \min\bigl( w+1, 2(e^w-1)\bigr) \leq3 w \|h'\|.
\end{eqnarray*}
This proves (\ref{exp-t1-1}) and hence Theorem \ref{exp-t1}.

\section*{Acknowledgments}
The authors thank Larry Goldstein for helping on the exponential
approximation and thank
an anonymous referee and an Associate Editor for their helpful comments.

%

%
\printaddresses

\end{document}